\documentclass[12pt,english,a4paper]{amsart}
\usepackage[T1]{fontenc}
\usepackage[latin1]{inputenc}
\usepackage{geometry}
\usepackage{amssymb}

\makeatletter

\let\SF@@footnote\footnote
\def\footnote{\ifx\protect\@typeset@protect
    \expandafter\SF@@footnote
  \else
    \expandafter\SF@gobble@opt
  \fi
}
\expandafter\def\csname SF@gobble@opt \endcsname{\@ifnextchar[%]
  \SF@gobble@twobracket
  \@gobble
}
\edef\SF@gobble@opt{\noexpand\protect
  \expandafter\noexpand\csname SF@gobble@opt \endcsname}
\def\SF@gobble@twobracket[#1]#2{}

 \theoremstyle{plain}    
 \newtheorem{thm}{Theorem}[section]
 \numberwithin{equation}{section} %% Comment out for sequentially-numbered
 \numberwithin{figure}{section} %% Comment out for sequentially-numbered
 \theoremstyle{plain}
 \theoremstyle{remark}    
 \newtheorem{acknowledgement}[thm]{Acknowledgement} 
 \theoremstyle{plain}    
 \newtheorem{lem}[thm]{Lemma} %%Delete [thm] to re-start numbering
 \theoremstyle{remark}
 \newtheorem{rem}[thm]{Remark}
 \theoremstyle{plain}    
 \newtheorem{cor}[thm]{Corollary} %%Delete [thm] to re-start numbering
 \theoremstyle{definition}
  \newtheorem{example}[thm]{Example}

\def\makebbb#1{
    \expandafter\gdef\csname#1\endcsname{
        \ensuremath{\Bbb{#1}}}
}
\makebbb{R}
\makebbb{N}
\makebbb{Z}
\makebbb{C}
\makebbb{G}
\makebbb{E}
\makebbb{H}
\makebbb{P}
\makebbb{B}
\makebbb{K}

\usepackage{babel}
\makeatother
\begin{document}

\title{Bergman kernels and equilibrium measures for ample line bundles }

\author{Robert Berman}

\email{robertb@math.chalmers.se}

\curraddr{Institut Fourier, 100 rue des Maths, BP 74, 38402 St Martin d'Heres
(France)}

\begin{abstract}
Let $L$ be an ample holomorphic line bundle over a compact complex
Hermitian manifold $X.$ Any fixed smooth hermitian metric $\phi$
on $L$ induces a Hilbert space structure on the space of global global
holomorphic sections with values in the $k$th tensor power of $L.$
In this paper various convergence results are obtained for the corresponding
Bergman kernels. The convergence is studied in the large $k$ limit
and is expressed in terms of the equilibrium metric associated to
the fixed metric $\phi,$ as well as in terms of the Monge-Ampere
measure of the metric $\phi$ itself on a certain support set. It
is also shown that the equilibrium metric has Lipschitz continuous
first derivatives. These results can be seen as generalizations of
well-known results concerning the case when the curvature of the fixed
metric $\phi$ is positive (the corresponding equilibrium metric is
then simply $\phi$ itself). 

\tableofcontents{}
\end{abstract}
\maketitle

\section{Introduction}

Let $L$ be an ample holomorphic line bundle over a compact complex
manifold $X$ of dimension $n.$ Fix an Hermitian fiber metric, denoted
by $\phi,$ on $L$ which is smooth and a smooth volume form on $\omega_{n}$
on $X.$ The curvature form of the metric $\phi$ may be written as
$dd^{c}\phi$ (see section \ref{sub:Notation-and-setup} for definitions
and further notation). Denote by $\mathcal{H}(X,L^{k})$ the Hilbert
space obtained by equipping the space $H^{0}(X,L^{k})$ of global
holomorphic sections with values in a tensor power $L^{k}$ with the
norm induced by the given smooth metric $\phi$ on $L$ and the volume
form $\omega_{n}.$ The \emph{Bergman kernel} of the Hilbert space
$\mathcal{H}(X,L^{k})$ is the integral kernel of the orthogonal projection
from the space of all smooth sections with values in $L^{k}$ onto
$\mathcal{H}(X,L^{k}).$ It may be represented by a holomorphic section
$K_{k}(x,y)$ of the pulled back line bundle $L^{k}\boxtimes\overline{L}^{k}$
over $X\times\overline{X}$ (formula \ref{def: K}).

In the case when the curvature form $dd^{c}\phi$ is \emph{positive,}
the asymptotic properties of the Bergman kernel $K_{k}(x,y)$ as $k$
tends to infinity has been studied thoroughly with numerous applications
in complex geometry and mathematical physics. For example, $K_{k}(x,y)$
admits a complete local asymptotic expansion in powers of $k;$ the
Tian-Zelditch-Catlin expansion (see \cite{ze,b-b-s} and references
therein). The point is that when the curvature form $dd^{c}\phi$
is globally positive, the Bergman kernel asymptotics at a fixed point
may be localized and hence only depend (up to negligable terms) on
the covariant derivatives of $dd^{c}\phi$ at the fixed point. 

The aim of the present paper is to study the case of a general smooth
metric $\phi$ on $L,$ where global effects become important and
where there appears to be very few previous general results. We will
consider three natural \emph{positive} measures on $X$ associated
to the setup introduced above. First the \emph{equilibrium measure}\[
(dd^{c}\phi_{e})^{n}/n!,\]
 where \emph{}$\phi_{e}$ is the \emph{equilibrium metric} defined
by the upper envelope \ref{eq:extem metric}, and then the large $k$
limit of the measures \begin{equation}
k^{-n}B_{k}\omega_{n},\label{eq:intro B meas}\end{equation}
 where $B_{k}(x):=K_{k}(x,x)e^{-k\phi}$ will be referred to as the
\emph{Bergman function} and of the measure \[
(dd^{c}(k^{-1}\textrm{ln\,$K_{k}(x,x)))^{n}/n!$, }\]
often referred to as the $k$th \emph{Bergman volume form} on $X$
associated to $(L,\phi).$

It is not hard to see that the total integrals of all three measures
coincide (and equal the total integral over $X$ of the (possibly
non-positive) form $(dd^{c}\phi)^{n}$). The main point of the present
paper is to show the corresponding \emph{local} statement. In fact,
all three measures will be shown to coincide with the measure \[
1_{D}(dd^{c}\phi)^{n}/n!\]
where $1_{D}$ is the characteristic function the $D$ in $X$ where
$\phi_{e}=\phi$ (corollary \ref{cor:equil meas}, theorem \ref{thm:B in L1}
and \ref{thm:ln K}). In the case when the metric $\phi$ has a semi-positive
curvature form, $\phi_{e}=\phi,$ i.e. the set $D$ equals all of
$X.$ The main results may be suggestively formulated in the following
form: \[
K_{k}(x,x):=(k^{n}\det(dd^{c}\phi_{e})(x)+...)e^{k\phi_{e}(x,x)},\]
 where the dots indicate terms of lower order in $k.$ On the set
$D,$ the \emph{}equilibrium metric \emph{$\phi_{e}$} may be replaced
by $\phi$ and away from the set $D$ the Bergman function $K_{k}(x,x)e^{-k\phi(x,x)}$
becomes exponentially small in $k$. Moreover, for any interior point
of $D$ where the curvature form of $\phi$ is positive we will show
(theorem \ref{thm:asymp expansion}) that the Bergman kernel $K_{k}(x,y),$
i.e. with two arguments, admits a complete local asymptotic expansion
in powers of $k,$ such that the coefficients of the corresponding
symbol expansion coincide with the Tian-Zelditch-Catlin expansion
for a \emph{positive} Hermitian holomorphic line bundle. Moreover,
it will be shown (theorem \ref{thm:k as meas}) that globally on $X\times X$
the following asymptotics hold for the point-wise norm $\left|K_{k}(x,y)\right|_{k\phi}:$\[
\begin{array}{lr}
k^{-n}\left|K_{k}(x,y)\right|_{k\phi}^{2}\omega_{n}(x)\wedge\omega_{n}(y)\rightarrow\Delta\wedge1_{D\cap X(0)}(dd^{c}\phi)^{n}/n!\end{array},\]
 weakly as measures on $X\times X$, where $\Delta$ is the current
of integration along the diagonal in $X\times X.$ A crucial step
in the present approach is to first show the $\mathcal{C}^{1,1}$-regularity
of the equilibrium metric $\phi_{e}$ (theorem \ref{thm:reg}), which
is of independent interest.

Finally, the setup above will be adapted to the case when the Hilbert
space $\mathcal{H}(X,L^{k})$ is replaced by the subspace of all sections
vanishing to high order along a fixed divisor in $X.$ The point is
that even in the case when the curvature form of the metric $\phi$
is positive, the introduction of the divisor is essentially equivalent
to studying a (singular) metric with \emph{negative} curvature concentrated
along the divisor.

\subsection{Comparison with previous results}

The present paper can be seen as a global geometric version of the
situation recently studied in \cite{berm4}, where the role of the
Hilbert space $\mathcal{H}(X,L^{k})$ was played by the space of all
polynomials in $\C^{n}$ of total degree less than $k,$ equipped
with a weighted norm. In fact, apart from the $\mathcal{C}^{1,1}$-regularity
of the equilibrium metric $\phi_{e},$ the adaptation to the present
setting is fairly straight forward. The proof of the $\mathcal{C}^{1,1}$-regularity
is partly modeled on the proof of Bedford-Taylor \cite{b-t,ko} for
$\mathcal{C}^{1,1}-$regularity of the solution of the Dirichlet problem
(with smooth boundary data) for the complex Monge-Ampere equation
in the unit-ball in $\C^{n}.$ The result should also be compared
to various $\mathcal{C}^{1,1}-$results for boundary value problems
for complex Monge-Ampere equations on manifolds with boundary \cite{ch,c-t},
intimately related to the study of the geometry of the space of Kähler
metrics on a Kähler manifold (see also \cite{p-s,bern2} for other
relations to Bergman kernels in the latter context). However, the
present situation rather corresponds to a \emph{free} boundary value
problem (compare remark \ref{rem:free}). 

In the case of sections vanishing along a fixed divisor (section \ref{sec:Sections-vanishing-along})
and under the further assumtion that the curvature form $dd^{c}\phi$
is positive, similar results have independently been obtained by Julien
Keller, Gabor Szekelyhidi and Richard Thomas \cite{kst}, but with
a different algebro-geometric characerization of the set $D_{Z}$
in formula \ref{eq:def of DZ}. 

As in the case studied in \cite{berm2} (part 1), where the curvature
form of the metric $\phi$ was assumed to be semi-positive the present
approach to the Bergman kernel asymptotics is based on the use of
{}``local holomorphic Morse-inequalities'', which are local version
of the global ones introduced by Demailly \cite{d1}. These inequalities
are then combined with some global pluripotential theory, based on
the recent work \cite{g-z} by Guedj-Zeriahi. Further references and
comments on the relation to the study of random polynomials (and holomorphic
sections), random eigenvalues of normal matrices and various diffusion-controlled
growth processes studied in the physics literature can be found in
\cite{berm4}.

\subsection{Further generalizations}

It can be further shown that the main results in this paper may be
generalized to \emph{any} line bundle $L$ over a Kähler manifold.
The results are then closely related to the study of the volume of
a line bundle over a Kähler manifold \cite{bo}. The main general
case is when the bundle $L$ is big (i.e. the dimension of $H^{0}(X,L^{k})$
is of the order $k^{n}$). Then all the results obtained in the ample
case still hold if the Monge-Ampere measure $(dd^{c}\phi)^{n}$ is
replaced by $1_{X-F}(dd^{c}(1_{X-F}\phi))^{n},$ where $F$ is a certain
analytic variety in $X,$ naturally associated to $L.$ The main (technical)
point is that the proof of the $\mathcal{C}^{1,1}-$ regularity of
$\phi_{e}$ goes through on the complement of $F.$ Finally, in the
case when $L$ is not big it can be shown that the convergence results
of the Bergman kernels simply say that $0=0.$ The details will appear
elsewhere.

\begin{acknowledgement}
It is a pleasure to thank Jean-Pierre Demailly and Sebastian Boucksom
for several illuminating discussions on the topic of the present paper.
In particular, the approach to the proof of the regularity theorem
\ref{thm:reg} was suggested by Jean-Pierre Demailly. Also thanks
to Julien Keller for informing me about the work \cite{kst}.
\end{acknowledgement}

\subsection{\label{sub:Notation-and-setup}General notation%
\footnote{general references for this section are the books \cite{gr-ha,de4}.%
}}

Let $(L,\phi)$ be an Hermitian holomorphic line bundle over a compact
complex manifold $X.$ The fixed Hermitian fiber metric on $L$ will
be denoted by $\phi.$ In practice, $\phi$ is considered as a collection
of local \emph{smooth} functions. Namely, let $s^{U}$ be a local
holomorphic trivializing section of $L$ over an open set $U$ then
locally, $\left|s^{U}(z)\right|_{\phi}^{2}=:e^{-\phi^{U}(z)},$ where
$\phi^{U}$ is in the class $\mathcal{C}^{2},$ i.e. it has continuous
derivatives of order two. If $\alpha_{k}$ is a holomorphic section
with values in $L^{k},$ then over $U$ it may be locally written
as $\alpha_{k}=f_{k}^{U}\cdot(s^{U})^{\otimes k},$ where $f_{k}^{U}$
is a local holomorphic function. In order to simplify the notation
we will usually omit the dependence on the set $U.$ The point-wise
norm of $\alpha_{k}$ may then be locally expressed as\begin{equation}
\left|\alpha_{k}\right|_{k\phi}^{2}=\left|f_{k}\right|^{2}e^{-k\phi}.\label{eq:ptwise norm}\end{equation}
 The canonical curvature two-form of $L$ is the global form on $X,$
locally expressed as $\partial\overline{\partial}\phi$ and the normalized
curvature form $i\partial\overline{\partial}\phi/2\pi=dd^{c}\phi$
(where $d^{c}:=i(-\partial+\overline{\partial})/4\pi)$ represents
the first Chern class $c_{1}(L)$ of $L$ in the second real de Rham
cohomology group of $X.$ The curvature form of a smooth metric is
said to be \emph{positive} at the point $x$ if the local Hermitian
matrix $(\frac{\partial^{2}\phi}{\partial z_{i}\partial\bar{z_{j}}})$
is positive definite at the point $x$ (i.e. $dd^{c}\phi_{x}>0).$
This means that the curvature is positive when $\phi(z)$ is strictly
\emph{plurisubharmonic} i.e. strictly subharmonic along complex lines.
We let \[
X(0):=\left\{ x\in X:\, dd^{c}\phi_{x}>0\right\} \]
 A line bundle $L$ is \emph{ample} if admits some smooth metric whose
curvature is positive on all of $X.$

More generally, a metric $\phi'$ on $L$ is called (possibly) \emph{singular}
if $\left|\phi'\right|$ is locally integrable. Then the curvature
is well-defined as a $(1,1)-$current on $X.$ The curvature current
of a singular metric is called \emph{positive} if $\phi'$ may be
locally represented by a plurisubharmonic function (in particular,
$\phi'$ takes values in $[-\infty,\infty[$ and is upper semi-continuous
(u.s.c)). In particular , any section $\alpha_{k}$ as above induces
such a singular metric on $L,$ locally represented by $\phi'=\frac{1}{k}\ln\left|f_{k}\right|^{2}.$
If $Y$ is a complex manifold we will denote by $PSH(Y)$ and $SPSH(Y)$
the space of all plurisubharmonic and strictly plurisubharmonic functions,
respectively. 

Fixing an Hermitian metric two-form $\omega$ on $X$ (with associated
volume form $\omega_{n})$ the Hilbert space $\mathcal{H}(X,L^{k})$
is defined as the space $H^{0}(X,L^{k})$ with the norm 

\begin{equation}
\left\Vert \alpha_{k}\right\Vert _{k\phi}^{2}(=\int_{X}\left|f_{k}\right|^{2}e^{-k\phi(z)}\omega_{n}),\label{eq:norm restr}\end{equation}
using a suggestive notation in the last equality (compare formula
\ref{eq:ptwise norm}).

\section{\label{sec:Equilibrium-measures-for}Equilibrium measures for line
bundles}

Let $L$ be a line bundle over a compact complex manifold $X.$ Given
a smooth metric $\phi$ on $L$ the corresponding {}``equilibrium
metric'' $\phi_{e}$ is defined as the envelope 

\begin{equation}
\phi_{e}(x)=\sup\left\{ \widetilde{\phi}(x):\,\widetilde{\phi}\in\mathcal{L}_{(X,L)},\,\widetilde{\phi}\leq\phi\,\,\textrm{on$\, X$}\right\} .\label{eq:extem metric}\end{equation}
 where $\mathcal{L}_{(X,L)}$ is the class consisting of all (possibly
singular) metrics on $L$ with positive curvature form. If $\phi_{e}$
is not u.s.c it should be replace it with its u.s.c regularization.
Then $\phi_{e}$ is also in the class $\mathcal{L}_{(X,L)}$ \cite{g-z}.
In the case when $\phi_{e}$ is locally bounded, the corresponding
\emph{equilibrium measure} is defined as the Monge-Ampere measure
$(dd^{c}\phi_{e})^{n}/n!$ (see \cite{g-z} for the definition of
the Monge-Ampere measure of a locally bounded metric, based on the
work \cite{b-t} in $\C^{n}$). When $L$ is ample $\phi_{e}$ is
clearly locally bounded and it can be shown directly, by adapting
the corresponding proof in $\C^{n}$ (see the appendix in \cite{s-t})
that $(dd^{c}\phi_{e})^{n}/n!$ vanishes on the complement of the
set. \begin{equation}
D:=\{\phi_{e}=\phi\}\subset X\label{eq:def of D}\end{equation}
However, the previous vanishing will also be a corollary of theorem
\ref{thm:B in L1} below.

\subsection{$\mathcal{C}^{1,1}-$regularity for ample line bundles}

In this section we will prove that the equilibrium metric $\phi_{e}$
associated to a smooth metric on an \emph{ample} line bundle $L$
is locally in the class $\mathcal{C}^{1,1}.$ As in \cite{berm4},
where the manifold $X$ was taken as $\C^{n},$ the proof is modeled
on the proof of Bedford-Taylor \cite{b-t,ko,de3} for $\mathcal{C}^{1,1}-$regularity
of the solution of the Dirichlet problem (with smooth boundary data)
for the complex Monge-Ampere equation in the unit-ball in $\C^{n}.$
However, as opposed to $\C^{n}$ and the unit-ball a generic compact
Kähler manifold $X$ has no global holomorphic vector fields. In order
to circumvent this difficulty we will reduce the regularity problem
on $X$ to a problem on the pseudoconvex manifold $Y,$ where $Y$
is the total space of the dual line bundle $L^{*},$ identifying the
base $X$ with its embedding as the zero-section in $Y.$ To any given
(possibly singular) metric $\phi$ on $L$ we may associate the {}``squared
norm function'' $h_{\phi}$ on $Y,$ where locally \[
h_{\phi}(z,w)=\left|w\right|^{2}\exp(\phi(z)).\]
In this way we obtain a bijection \begin{equation}
\mathcal{L}_{(X,L)}\leftrightarrow\mathcal{L}_{Y},\,\,\,\phi\mapsto h_{\phi}\label{eq:bijection}\end{equation}
where $\mathcal{L}_{Y}$ is the class of all positively 2-homogeneous
plurisubharmonic functions on $Y:$ \begin{equation}
\mathcal{L}_{Y}:=\{ h\in PSH(Y):\, h(\lambda\cdot)=\left|\lambda\right|^{2}h(\cdot)\},\label{eq:class ly}\end{equation}
using the natural multiplicative action of $\C^{*}$ on the fibers
of $Y$ over $X.$ Now we define \begin{equation}
h_{e}:=\sup\left\{ h\in\mathcal{L}_{Y}:\, h\leq h_{\phi}\,\,\textrm{on$\, X$}\right\} ).\label{def: h e}\end{equation}
 Then clearly, $h_{e}$ corresponds to the equilibrium metric $\phi_{e}$
under the bijection \ref{eq:bijection}. The following lemma will
allow us to {}``homogenize'' plurisubharmonic functions. 

\begin{lem}
\label{lem:homo}Suppose that the function $f$ is spsh, $f=0$ on
$X,$ and $f$ is $S^{1}-$invariant in a neighbourhood of the sublevel
set $f^{-1}[-\infty,c],$ $c>0.$ Then there is a function $\widetilde{f}$
in the class $\mathcal{L}_{Y}$ such that $\widetilde{f}=f$ on the
level set $f^{-1}(c).$ 
\end{lem}
\begin{proof}
First observe that $f^{-1}(c)$ is an $S^{1}-$bundle subbundle of
$Y$ over $X.$ Indeed, since $f$ is assumed to be a spsh function
on $f^{-1}[-\infty,c]$ it follows, from the maximum principle applied
to discs in each fiber, that $f$ is strictly increasing along the
fibers of $Y$ over $X.$ Hence, since $f=0$ on the zero-section
and $c>0$ any $y$ in $Y-X$ may be written in a unique way as $y=r\sigma,$
where $f(\sigma)=c$ and $r>0$ and we may define $\widetilde{f}$
by \[
\widetilde{f}(r\sigma):=r^{2}\widetilde{f}(\sigma).\]
Finally, to see that $\widetilde{f}$ is psh note that since $f$
is strictly increasing along the fibers we have that $\widetilde{f}-c$
is a defining function for the domain $f^{-1}[-\infty,c],$ which
is pseudoconvex, since $f$ is psh. In the case when $f$ is smooth
and $c$ is a regular value it follows that $(dd^{c}\widetilde{f})\geq0$
along the holomorphic subbundle $T^{1,0}(f^{-1}[-\infty,c]).$By homogeneity
this means that $\widetilde{f}$ is psh on all of $Y.$ Finally, the
general case may be obtained by a local approximation argument.
\end{proof}
The next (essentially well-known) lemma provides the vector fields
needed in the approach of Bedford-Taylor:

\begin{lem}
\label{lem:exist of vf}Assume that the line bundle $L$ is ample.
For any given point $y_{0}$ in $Y-X$ there are global holomorphic
vector fields $V_{1},...V_{n+1}$ (i.e. elements of $H^{0}(Y,TY)$)
such that their restriction to $y_{0}$ span the tangent space $TY_{y_{0}}$
and such that $V_{i}$ vanish on $X.$ 
\end{lem}
\begin{proof}
It is well-known \cite{de4} that on a Stein-space $M$ any holomorphic
coherent sheaf is globally generated (i.e. it has the spanning property
stated in the lemma). Hence the lemma could be obtained by observing
that $Y$ may be blown-down to a Stein space after contracting $X$
to a point so that the push-forward of $TY$ becomes a coherent sheaf.
But for completeness we give a somewhat more explicit argument. First
note that $Y$ may be compactified by the following fiber-wise projectivized
vector bundle:\[
\widehat{Y}:=\P(L^{*}\oplus\underline{\C}),\]
 where $\underline{\C}$ denotes the trivial line bundle over $X.$
Denote by $\mathcal{O}(1)$ the line bundle over $\widehat{Y}$ whose
restriction to each fiber (i.e. a one-dimensional complex space $\P^{1})$
is the induced hyperplane line bundle. Next, observe that the line
bundle \[
\widehat{L}:=(\pi^{*}(L)\otimes\mathcal{O}(1))\]
 over $\widehat{Y}$ is ample if $L$ is, where $\pi$ denotes the
natural projection from $\widehat{Y}$ to $X.$ Indeed, any given
smooth metric $\phi_{+}$ on $L$ with positive curvature induces
a metric $\widehat{\phi}$ on $\widehat{L}$ with positive curvature,
expressed as \[
\widehat{\phi}=\pi^{*}\phi_{+}+\ln(1+h_{\phi_{+}})\]
on $Y$ (extending to $\widehat{Y}$). Now it is well-known (for example
using Hörmander's $L^{2}-$estimates \cite{de4}), that for any ample
line bundle $\widehat{L}$ and holomorphic vector bundle $E$ on a
compact manifold $\widehat{Y}$ the bundle $E\otimes\widehat{L}^{k_{0}}$is
globally generated for $k_{0}$ sufficiently large. Setting $E=T\widehat{Y}$
and restricting to $Y$ in $\widehat{Y}$ shows that $TY\otimes\pi^{*}(L){}^{k_{0}}$
is globally generated on $Y$ (since $\mathcal{O}(1)$ is trivial
on $Y).$ Finally, observe that \[
\pi^{*}(L)=(\pi^{*}(L^{*}))^{-1}=[X]^{-1},\]
 where $[X]$ is the divisor in $Y$ determined by the embedding of
$X$ as the base. Indeed, $X$ is embedded as the zero-set of the
tautological section of $\pi^{*}(L^{*})$ over $Y(=L^{*}).$ Hence,
the sections of $TY\otimes\pi^{*}(L){}^{k_{0}}$ may be identified
with sections in $TY$ vanishing to order $k_{0}$ on $X.$ This finishes
the proof of the lemma.
\end{proof}
For any given smooth vector field $V$ on $Y$ and compact subset
$K$ of $Y,$ we denote by $exp(tV)$ the corresponding flow which
is well-defined for any {}``time'' $t$ in $[0,t_{K}],$ i.e. the
family of smooth maps indexed by $t$ such that \begin{equation}
\frac{d}{dt}f(exp(tV)(y))=df[V]_{exp(tV)(y)}\label{eq:def of flow}\end{equation}
for any smooth function $f$ and point $y$ on $Y.$ We will also
use the notation $exp(V):=exp(1V).$ In the following we will just
be interested in an arbitrarily small neighbourhood of $X$ in $Y$
and we will tacitly take $t$ sufficiently small in order that the
flow exists (or equivalently, rescale $V).$ 

Now fix a point $y_{0}$ in $Y-X.$ Combining the previous lemma with
the inverse function theorem gives local {}``exponential'' holomorphic
coordinates centered at $y_{0},$ i.e a local biholomorphism \[
\C^{n+1}\rightarrow U(y),\,\,\,\lambda\mapsto exp((V(\lambda)(y_{0}),\,\,\, V(\lambda):=\sum\lambda_{i}V_{i})\]
Using that the vector fields $V_{i}$ necessarily also span $TY_{y_{1}}$
for $y_{1}$ close to $y_{0}$ it can be checked that in order to
prove that a function $f$ is locally Lipschitz continuous on a compact
subset of $Y$ it is enough to, for each fixed point $y_{0},$ prove
an estimate of the form \begin{equation}
\left|f(exp(V(\lambda))(y_{0}))-f(y_{0})\right|\leq C\left|\lambda\right|\label{eq:cond for lip}\end{equation}
for some constant $C$ only depending on the function $f.$ 

\begin{thm}
\label{thm:reg}Suppose that $L$ is an ample line bundle and that
the given metric $\phi$ on $L$ is smooth (i.e. in the class $\mathcal{C}^{2}).$
Then 

(a) $\phi_{e}$ is locally in the class $\mathcal{C}^{1,1},$ i.e.
$\phi_{e}$ is differentiable and all of its first partial derivatives
are locally Lipschitz continuous. 

(b) The Monge-Ampere measure of $\phi_{e}$ is absolutely continuous
with respect to any given volume form and coincides with the corresponding
$L_{loc}^{\infty}$ $(n,n)-$form obtained by a point-wise calculation:
\begin{equation}
(dd^{c}\phi_{e})^{n}=\det(dd^{c}\phi_{e})\omega_{n}\label{eq:ptwise repr of equil meas}\end{equation}

(c) the following identity holds almost everywhere on the set $D=\{\phi_{e}=\phi\}:$
\begin{equation}
\det(dd^{c}\phi_{e})=\det(dd^{c}\phi)\label{eq:monge on D}\end{equation}

\end{thm}
\begin{proof}
To prove $(a)$ it is, by the bijection \ref{eq:bijection}, equivalent
to prove that $h_{e}$ (defined by \ref{def: h e}) is locally $\mathcal{C}^{1,1}$
on $Y-X.$ Moreover, by homogeneity it is enough to show that there
is a neighbourhood $U$ of $X$ in $Y$ such that $h_{e}$ is locally
$\mathcal{C}^{1,1}$ on $U-X$ 

\emph{Step1:} \textbf{\emph{}}\emph{$h_{e}$ is locally Lipschitz
continuous on $Y-X.$ }

To see this fix a point $y_{0}$ in $Y-X$. For simplicity we first
assume that $h_{e}$ is \emph{strictly} plurisubharmonic on $Y-X$
(the assumption will be removed in the end of the argument). Let \[
g(y):=h_{e}(exp(V(\lambda))(y)),\,\,\,\widehat{g}(y):=\textrm{u.s.c}(\sup_{\theta\in[0,2\pi}g(e^{i\theta}y))\]
using the natural multiplicative action of $\C^{*}$ on the fibers
of $Y$ over $X$ and where u.s.c. denotes the upper-semicontinous
regularization. Then both $g$ and $\widehat{g}$ are psh functions
(using that the family $g(e^{i\theta}\cdot)$ of psh functions is
locally bounded). Note that\[
h_{e}(exp(V(\lambda))(y_{0})=:g(y_{0})\leq\widehat{g}(y_{0})=\widetilde{\widehat{g}}(y_{0}),\]
 where $\widetilde{\widehat{g}}$ is the function in the class $\mathcal{L}_{Y}$
obtained from lemma \ref{lem:homo} applied to $f=\widehat{g}$ and
$c=\widehat{g}(y_{0})$ (note that $\widehat{g}$ is still \emph{strictly}
plurisubharmonic, since the flow $exp(V(\lambda))$ fixes $X).$ Moreover,
since by definition $h_{e}\leq h_{\phi}$ we have the following bound
on the levelset $\widetilde{\widehat{g}}^{-1}(c):$ \begin{equation}
\widetilde{\widehat{g}}((y)\leq\sup_{\theta\in[0,2\pi]}h_{\phi}((exp(V(\lambda))(e^{i\theta}y)\leq\sup_{\theta\in[0,2\pi]}h_{\phi}(e^{i\theta}y)+C\left|\lambda\right|,\label{eq:pf thm reg}\end{equation}
using that $h_{\phi}$ is locally Lipschitz in the last inequality.
Indeed, since $h_{\phi}$ is in the class $\mathcal{C}^{1}$ the property
\ref{eq:def of flow} of the flow gives\[
h_{\phi}((exp(V(\lambda))(\zeta)-h_{\phi}(\zeta)=\int_{0}^{1}dh_{\phi}[V(\lambda)]_{exp(tV(\lambda)(y)}dt\]
Hence, since $V(\lambda):=\left|\lambda\right|(\sum\frac{\lambda_{i}}{\left|\lambda\right|}V_{i}),$
the constant $C$ in \ref{eq:pf thm reg} may be taken to be \[
C=\sup_{y\in K,i=1,..N+1}\left|dh_{\phi}[V_{i}]_{y}\right|\]
 for some compact neighbourhood $K$ of $X$ in $Y.$ Since $h_{\phi}$
is $S^{1}-$invariant and $\widetilde{\widehat{g}}^{-1}(c)$ is compact,
\ref{eq:pf thm reg} gives that \begin{equation}
\widetilde{\widehat{g}}-C\left|\lambda\right|\leq h_{\phi}\label{pf thm reg: bound on cand}\end{equation}
 on $\widetilde{\widehat{g}}^{-1}(c)$ and hence, by homogeneity,
on all of $Y.$ This shows that the function $\widetilde{\widehat{g}}-C\left|\lambda\right|$
is a contender for the supremum in the definition \ref{def: h e}
of $h_{e}$ and hence bounded by $h_{e}.$ All in all we get that
\[
h_{e}(exp(V(\lambda))(y_{0})\leq\widetilde{\widehat{g}}(y_{0})\leq h_{\phi}(y_{0})+C\left|\lambda\right|.\]
The other side of the inequality \ref{eq:cond for lip} for $f=h_{e}$
is obtained after replacing $\lambda$ by $-\lambda.$ Finally, in
order to remove the simplifying assumption that $h_{e}$ be \emph{strictly}
plurisubharmonic on $Y-X$ we apply the previous argument to get the
same bounds on \[
h_{\delta}:=h_{e}(1-\delta)+\delta h_{+}\]
where $h_{+}$ is spsh and in the class $\mathcal{L}_{Y}$ (if $L$
is ample than $h_{+}$ clearly exists). Finally, letting $\delta$
tend to zero, finishes the proof of Step 1. 

\emph{Step2: $dh_{e}$ exists and is locally Lipschitz continuous
on $Y-X.$}

Following the exposition in \cite{de3} it is enough to prove the
following inequality:

\begin{equation}
h_{e}(exp(V(\lambda))(y_{0})+h_{e}(exp(V(-\lambda))(y_{0})-2h_{e}(y_{0})\leq C\left|\lambda\right|^{2},\label{eq:pf of thm reg ineq 2}\end{equation}
where the constant only depends on the second derivatives of $h_{\phi}$
on some compact subset of $Y.$ Indeed, given this inequality (combined
with the fact that $h_{e}$ is psh) a Taylor expansion of degree $2$
gives the following bound close to $y_{0}$ for a smooth approximation
$h_{\epsilon}$ of $h_{e}:$ \[
\left|D^{2}h_{\epsilon}\right|\leq C\]
 where $h_{\epsilon}:=h_{e}*\chi_{\epsilon},$ using a a local regularizing
kernel $\chi_{\epsilon}$ and where $D^{2}h_{\epsilon}$ denotes the
real Hessian matrix of $h_{\epsilon}.$ Letting $\epsilon$ tend to
$0$ then proves Step 2. Finally, to see that the inequality \ref{eq:pf of thm reg ineq 2}
holds we apply the argument in Step 1 after replacing $g$ by the
psh function \[
f(y):=(h_{e}(exp(V(\lambda))(y)+h_{e}(exp(V(-\lambda))(y))/2\]
 to get \[
f(y)\leq\widetilde{\widehat{f}}((y)\leq\sup_{\theta\in[0,2\pi]}(h_{\phi}((exp(V(\lambda))(e^{i\theta}y)+(exp(V(-\lambda))(e^{i\theta}y))/2\]
Next, observe that for each fixed $\theta$ the function $h_{\phi}(e^{i\theta}y)$
is in the class $\mathcal{C}^{2}.$ Hence, a Taylor expansion of degree
$2$ gives \[
\widetilde{\widehat{f}}((y)\leq\sup_{\theta\in[0,2\pi]}((h_{\phi}(e^{i\theta}y))+C\left|\lambda\right|^{2})=h_{\phi}(y)+C\left|\lambda\right|^{2})\]
 where the constant $C$ may be taken as a constant times $\sup_{K}\left|D^{2}h_{\phi}\right|.$
This shows that $\widetilde{\widehat{f}}-C\left|\lambda\right|^{2}$
is a contender for the supremum in the definition \ref{def: h e}
of $h_{e}$ and hence bounded by $h_{e}.$ All in all we obtain that
\[
f(y_{0})\leq h_{\phi}(y_{0})+C\left|\lambda\right|^{2},\]
 which proves the inequality \ref{eq:pf of thm reg ineq 2}, finishing
the proof of Step2.

(b) By the $\mathcal{C}^{1,1}-$regularity, the derivatives $\frac{\partial^{2}\phi}{\partial z_{i}\partial\bar{z_{j}}}\phi_{e}$
are in $L_{loc}^{\infty}$ and it is well-known that this implies
the identity \ref{eq:ptwise repr of equil meas} for the Monge-Ampere
measure. Finally, to see that \ref{eq:monge on D} holds, it is enough
to prove that locally\[
\frac{\partial^{2}\phi}{\partial z_{i}\partial\bar{z_{j}}}(\phi_{e}-\phi)=0\]
almost everywhere on $D=\{\phi_{e}=\phi\}.$ To this end we apply
a calculus lemma in \cite{k-s} (page 53) to the $\mathcal{C}^{1,1}-$function
$\phi_{e}-\phi$ (following the approach in \cite{h-m}), which even
gives the corresponding identity between all real second order partial
derivatives almost everywhere on $D.$
\end{proof}
\begin{rem}
\label{rem:free}Fix a metric $\phi_{+}$ on $L$ with positive curvature.
Then $\omega_{+}:=dd^{c}\phi_{+}$ is a Kähler metric on $X$ and
the fixed metric $\phi$ on $L$ may be written as $\phi=u+\phi_{+},$
where $u$ is a smooth \emph{function} on $X.$ Now the pair $(u_{e},M)$
where $u_{e}:=\phi_{e}-\phi_{+}$ and $M$ is the set $X-D,$ may
be interpreted as a {}``weak'' solution to the following \emph{free}
boundary value problem of Monge-Ampere type%
\footnote{since there is a priori no control on the regularity of the set $M,$
it does not really make sense to write $\textrm{$\partial M$}$ and
the boundary condition should hence be interpreted in a suitable {}``weak''
sense.%
}: \[
\begin{array}{rclr}
(dd^{c}u_{e}+\omega_{+})^{n} & = & 0 & \textrm{on\,$M$}\\
u_{e} & = & u & \textrm{on\,$\partial M$}\\
du_{e} & = & du\end{array}\]
The point is that, since the equations are overdetermined, the set
$M$ is itself part of the solution. In \cite{h-m} the $\mathcal{C}^{1,1}-$regularity
of $\phi_{e}$ in the case when $X=\C$ (corresponding to the setup
in \cite{s-t}) was deduced from the regularity of a free boundary
value problem. 
\end{rem}

\section{\label{sec:Bergman-kernels}Bergman kernel asymptotics}

Denote by $\mathcal{H}(X,L^{k})$ the Hilbert space obtained by equipping
the vector space $H^{0}(X,L^{k})$ with the norm \ref{eq:norm restr}
induced by the given smooth metric $\phi$ on $L$ and the volume
form $\omega_{n}.$ Let $(\psi_{i})$ be an orthonormal base for $\mathcal{H}(X,L^{k}).$
The \emph{Bergman kernel} of the Hilbert space $\mathcal{H}(X,L^{k})$
is the integral kernel of the orthogonal projection from the space
of all smooth sections with values in $L^{k}$ onto $\mathcal{H}(X,L^{k}).$
It may be represented by the holomorphic section \begin{equation}
K_{k}(x,y)=\sum_{i}\psi_{i}(x)\otimes\overline{\psi_{i}(y)}.\label{def: K}\end{equation}
 of the pulled back line bundle $L^{k}\boxtimes\overline{L}^{k}$
over $X\times\overline{X}.$ The restriction of $K_{k}$ to the diagonal
is a section of $L^{k}\otimes\overline{L}^{k}$ and we let $B_{k}(x)=\left|K_{k}(x,x)\right|_{k\phi}(=\left|K_{k}(x,x)\right|e^{-k\phi(x)})$
be its point wise norm: \begin{equation}
B_{k}(x)=\sum_{i}\left|\psi_{i}(x)\right|_{k\phi}^{2}.\label{eq:def of B}\end{equation}
We will refer to $B_{k}(x)$ as the \emph{Bergman function} of $\mathcal{H}(X,L^{k}).$
It has the following extremal property:\begin{equation}
B_{k}(x)=\sup\left\{ \left|\alpha_{k}(x)\right|_{k\phi}^{2}:\,\,\alpha_{k}\in\mathcal{H}(X,L^{k}),\,\left\Vert \alpha_{k}\right\Vert _{k\phi}^{2}\leq1\right\} \label{(I)extremal prop of B}\end{equation}
 Moreover, integrating \ref{eq:def of B} shows that $B_{k}$ is a
{}``dimensional density'' of the space $\mathcal{H}(X,L^{k}):$
\begin{equation}
\int_{X}B_{k}\omega_{n}=\dim\mathcal{H}(X,L^{k})\label{eq:dim formel for B}\end{equation}
Now if $L$ is an ample line bundle, the dimension of $\mathcal{H}(X,L^{k})$
is explicitly known to the leading order in $k$ \cite{gr-ha} giving
\begin{equation}
\lim_{k\rightarrow\infty}\int_{X}k^{-n}B_{k}\omega_{n}=\int_{X}c_{1}(L)^{n}(=\int_{X}(dd^{c}\phi_{e})^{n}),\label{eq:int B as curv}\end{equation}
 where we have represented the first Chern class $c_{1}(L)$ by the
\emph{positive} curvature current $dd^{c}\phi_{e}$ of the equilibrium
metric $\phi_{e},$ using that $(dd^{c}\phi_{e})^{n}$ is well-defined
when $L$ is ample.

The following {}``local Morse inequality'' estimates $B_{k}$ point-wise
from above for a general bundle:

\begin{lem}
\emph{\label{lem:(Local-Morse-inequalities)}(Local Morse inequalities)}
Let $\phi$ be a smooth metric on a holomorphic line bundle $L$ over
a compact manifold $X.$ Then the following upper bound holds on $X:$
\[
k^{-n}B_{k}\leq C_{k}1_{X(0)}\det(dd^{c}\phi),\]
 where the sequence $C_{k}$ of positive numbers tends to one and
$X(0)$ is the set where $dd^{c}\phi>0.$
\end{lem}
See \cite{berm1} for the more general corresponding result for $\overline{\partial}-$harmonic
$(0,q)$- forms with values in a high power of an Hermitian line bundle.
The present case (i.e. $q=0)$ is a simple consequence of the mean-value
property of holomorphic functions applied to a poly-disc $\Delta_{k}$
of radius $\ln k/\sqrt{k}$ centered at the origin in $\C^{n}$ (see
the proof in \cite{berm2}). In fact, the proof gives the following
stronger local statement: \begin{equation}
\limsup_{k}k^{-n}\left|f_{k}(0)\right|^{2}e^{-k\phi(z)}/\left\Vert f_{k}\right\Vert _{k\phi,\Delta_{k}(z)}^{2}\leq1_{X(0)}(0)\det(dd^{c}\phi),\label{eq:local morse on ball}\end{equation}
 where $f_{k}$ is holomorphic function defined in a fixed neighbourhood
of the origin in $\C^{n}.$

The estimate in the previous lemma can be considerably sharpened on
the complement of $D$ (formula \ref{eq:def of D}), as shown by the
following lemma:

\begin{lem}
\label{lem:exponentia decay}Let $\phi$ be a smooth metric on a holomorphic
line bundle $L$ over a compact manifold $X.$ Then the following
inequality holds on all of $X$: \begin{equation}
B_{k}k^{-n}\leq C_{k}e^{-k(\phi-\phi_{e})}\label{eq:exp decay}\end{equation}
 where the sequence $C_{k}$ of positive numbers tends to $\sup_{X}\det(dd^{c}\phi).$
In particular, \begin{equation}
\lim\int_{D^{c}}k^{-n}B_{k}\omega_{n}=0\label{pf of thm B: int B on D compl}\end{equation}

\end{lem}
\begin{proof}
By the extremal property \ref{(I)extremal prop of B} of $B_{k}$
it is enough to prove the lemma with $B_{k}k^{-n}$replaced by $\left|\alpha_{k}\right|_{k\phi}^{2},$
locally represented by $\left|f_{k}\right|e^{-k\phi},$ for any element
$\alpha_{k}$ in $\mathcal{H}(X,L^{k})$ with global norm equal to
$k^{-n}.$ The Morse inequalities in the previous lemma give that\[
\left|f_{k}\right|^{2}e^{-k\phi}\leq C_{k}\]
 with $C_{k}$ as in the statement of the present lemma. Equivalently,
\[
\frac{1}{k}\ln\left|f_{k}\right|^{2}-\frac{1}{k}C_{k}\leq\phi\]
 Hence, the singular metric on $L$ determined by $\frac{1}{k}\ln\left|f_{k}\right|^{2}-\frac{1}{k}C_{k}$
is a candidate for the sup in the definition \ref{eq:extem metric}
of $\phi_{e}$ and is hence bounded by $\phi_{e}.$ Thus, \[
B_{k}k^{-n}=\left|f_{k}\right|^{2}e^{-k\phi}\leq C_{k}e^{k\phi_{e}}e^{-k\phi}.\]

\end{proof}
Finally, the vanishing \ref{pf of thm B: int B on D compl} follows
from the dominated convergence theorem, since the right hand side
in the previous inequality tends to zero precisely on the complement
of $D.$

\begin{thm}
\label{thm:B in L1}Let $L$ be an ample line bundle over $X$ and
let $B_{k}$ be the Bergman function of the Hilbert space $\mathcal{H}(X,L^{k}).$
Then\begin{equation}
k^{-n}B_{k}\rightarrow1_{D\cap X(0)}\det(dd^{c}\phi),\label{eq:l1 conv of B}\end{equation}
 in $L^{1}(X,\omega_{n}),$ where $X(0)$ is the set where $dd^{c}\phi>0$
and $D$ is the set \ref{eq:def of D}. 
\end{thm}
\begin{proof}
First observe that, by the exponential decay in lemma \ref{lem:exponentia decay},
\[
\lim_{k\rightarrow\infty}k^{-n}B_{k}(x)=0,\,\, x\in D^{c}\]
Next, observe that it is enough to prove that \begin{equation}
\lim_{k\rightarrow\infty}\int_{D}k^{-n}B_{k}\omega_{n}=\int_{D\cap X(0)}(dd^{c}\phi)^{n}/n!\label{pf of thm B: claim}\end{equation}
Indeed, given this equality the local Morse inequalities (lemma \ref{lem:(Local-Morse-inequalities)}),
then force the convergence \ref{eq:l1 conv of B} on the compact set
$D.$ The proof proceeds precisely as in \cite{berm2} (part 1, section
2).

Finally, to prove that \ref{pf of thm B: claim} does hold, first
note that 

\begin{equation}
\int_{X}(dd^{c}\phi_{e})^{n}=\lim_{k\rightarrow\infty}\int_{D}k^{-n}B_{k}\omega_{n}\leq\int_{D\cap X(0)}(dd^{c}\phi)^{n}/n!\label{eq:pf of thm B many ineq}\end{equation}
where we have combined formula \ref{eq:int B as curv} and \ref{pf of thm B: int B on D compl}
to get the equality and then used the local Morse inequalities (lemma
\ref{lem:(Local-Morse-inequalities)}) in the inequality (also using
the dominated convergence theorem). Finally, by formula \ref{eq:monge on D}
in theorem \ref{thm:reg} we may replace $\phi$ with $\phi_{e}$
in the right hand side, showing the the right hand side in the previous
equality is equal to $\int_{D\cap X(0)}(dd^{c}\phi_{e})^{n}.$ But
since $(dd^{c}\phi_{e})^{n}$ is a positive measure this can only
happen if all inequalities in \ref{eq:pf of thm B many ineq} are
actually equalities, which proves \ref{pf of thm B: claim} and finishes
the proof of the theorem. 
\end{proof}
Combining the previous theorem with the regularity theorem \ref{thm:reg}
now gives the following

\begin{cor}
\label{cor:equil meas}The equilibrium measure corresponding to the
smooth metric $\phi$ on $L$ is given by \[
(dd^{c}\phi_{e})^{n}/n!=1_{D\cap X(0)}(dd^{c}\phi)^{n}/n!,\]
where $D=\{\phi_{e}=\phi\}.$
\end{cor}

\subsection{The Bergman metric}

The Hilbert space $\mathcal{H}(X,L^{k})$ induces a metric on the
line bundle $L$ in the class $\mathcal{L}_{(X,L)}$ which may be
expressed as \[
k^{-1}\textrm{ln\,$K_{k}(x,x),$}\]
 often referred to as the $k$th Bergman metric on $L.$ If $L$ is
an ample line bundle, then this is the smooth metric on $L$ obtained
as the pull-back of the Fubini-Study metric on the hyperplane line
bundle $\mathcal{O}(1)$ over $\P^{N}(=\P\mathcal{H}(X,L^{k}))$ (compare
example \ref{exa:Pn} in section \ref{sec:Examples}) under the Kodaira
map \[
X\rightarrow\P\mathcal{H}(X,L^{k}),\,\,\, y\mapsto(\textrm{$\Psi_{1}(x):\Psi_{2}(x)...:\Psi_{N}(x))$ , }\]
for $k$ sufficiently large, where $\textrm{$(\Psi_{i})$}$ is an
orthonormal base for $\mathcal{H}(X,L^{k})$ \cite{gr-ha}. 

We will make use of the following well-known extension lemma, which
follows from the Ohsawa-Takegoshi theorem (compare \cite{bern2}):

\begin{lem}
\label{lem:o t}Let $F$ be a line bundle with a (possibly singular)
metric $\phi_{L}$ such that the curvature $dd^{c}\phi_{F}$ is positive
in the sense of currents and let $A$ be an ample line bundle. Then,
after possibly replacing $A$ by a sufficiently large tensor power,
the following holds: for any point $x$ in $X$ where $\phi\neq\infty,$
there is an element $\alpha$ in $H^{0}(X,F\otimes A)$ such that
\begin{equation}
\left|\alpha(x)\right|_{\phi_{F}+\phi_{A}}=1,\,\,\,\left\Vert \alpha\right\Vert _{X,\phi_{F}+\phi_{A}}\leq C.\label{eq:lem o t}\end{equation}
The constant $C$ is independent of the line bundle $F$ and the point
$x$ and depends only on a fixed smooth metric $\phi_{A}$ on $A.$ 
\end{lem}
Now we can prove the following theorem:

\begin{thm}
\label{thm:ln K}Let $L$ be an ample line bundle and let $K_{k}$
be the Bergman kernel of the Hilbert space $\mathcal{H}(X,L^{k}).$
Then the following convergence of Bergman metrics holds:\begin{equation}
k^{-1}\textrm{ln\,$K_{k}(x,x)\rightarrow\phi_{e}(x)$}\label{eq:conv of ln K in theorem ln K}\end{equation}
uniformly on $X$ (the rate of convergence is of the order $\ln k/k).$
In particular, the corresponding {}``Bergman volume forms'' converge
to the equilibrium measure: 
\end{thm}
\begin{equation}
(dd^{c}(k^{-1}\textrm{ln\,$K_{k}(x,x)))^{n}\rightarrow(dd^{c}\phi_{e})^{n}$}\label{eq:cong of monge in thm ln k}\end{equation}
weakly as measures.

\begin{proof}
In the following proof it will be convenient to let $C$ denote a
sufficiently large constant (which may hence vary from line to line).
First observe that taking the logarithm of the inequality \ref{eq:exp decay}
in lemma \ref{lem:exponentia decay} immediately gives the upper bound
\[
k^{-1}\textrm{ln\,$K_{k}(x,x)\leq\phi_{e}(x)+C\ln k/k$}\]
To get a lower bound, fix a point $x_{0}$ in $X$ and note that by
the extremal property \ref{(I)extremal prop of B} it is enough to
find a section $\alpha_{k}$ such that \begin{equation}
\left|\alpha_{k}(x_{0})\right|_{k\phi_{e}}\geq1/C,\,\,\,\left\Vert \alpha_{k}\right\Vert _{X,k\phi}\leq C.\label{pf of theorem ln K}\end{equation}
To this end we take the section $\alpha_{k}$ furnished by lemma the
previous lemma applied to $(F,\phi_{F})=(L^{k-k_{0}},(k-k_{0})\phi_{e})$
and $A=L^{k_{0}}$, for $k_{0}$ sufficiently large. Then $L^{k}=F\otimes A$
gets an induced metric \begin{equation}
\psi_{k}=(k-k_{0})\phi_{e}+\phi_{A}\label{pf of thm lnk: metric}\end{equation}
 and \[
\left|\alpha_{k}(y)\right|_{k\phi_{e}}\geq1/C,\,\,\,\left\Vert \alpha_{k}\right\Vert _{X,\psi_{k}}\leq C,\]
 where we have fixed a smooth metric $\phi_{A}$ on $L^{k_{0}}$ with
positive curvature such that $\phi_{A}\leq k_{0}\phi.$ Since, by
definition $\phi_{e}\leq\phi,$ this proves \ref{pf of theorem ln K}
and hence the theorem. 

The Monge-Ampere convergence \ref{eq:cong of monge in thm ln k} now
follows from the uniform convergence \ref{eq:conv of ln K in theorem ln K}
(see \cite{g-z}).
\end{proof}

\subsection{The full Bergman kernel }

Combining the convergence in theorem \ref{thm:B in L1} with the local
inequalities \ref{eq:local morse on ball}, gives the following convergence
for the point-wise norm of the full Bergman kernel $K_{k}(x,y).$
The proof is completely analogous to the proof of theorem 2.4 in part
1 of \cite{berm2}. 

\begin{thm}
\label{thm:k as meas}Let $L$ be an ample line bundle and let $K_{k}$
be the Bergman kernel of the Hilbert space $\mathcal{H}(X,L^{k}).$
Then \[
\begin{array}{lr}
k^{-n}\left|K_{k}(x,y)\right|_{k\phi}^{2}\omega_{n}(x)\wedge\omega_{n}(y)\rightarrow\Delta\wedge1_{D\cap X(0)}(dd^{c}\phi)^{n}/n!\end{array},\]
as measures on $X\times X$, in the weak {*}-topology, where $\Delta$
is the current of integration along the diagonal in $X\times X.$
\end{thm}
Finally, we will show that around any interior point of the set $D\cap X(0)$
the Bergman kernel $K_{k}(x,y)$ admits a complete local asymptotic
expansion in powers of $k,$ such that the coefficients of the corresponding
symbol expansion coincide with the Tian-Zelditch-Catlin expansion
(concerning the case when the curvature form of $\phi$ is positive
on all of $X;$ see \cite{b-b-s} and the references therein for the
precise meaning of the asymptotic expansion). We will use the notation
$\phi(x,y)$ for a fixed almost holomorphic-anti-holomorphic extension
of a local representation of the metric $\phi$ from the diagonal
$\Delta$ in $\C^{n}\times\C^{n},$ i.e. an extension such that the
anti-holomorphic derivatives in $x$ and the holomorphic derivatives
in $y$ vanish to infinite order on $\Delta.$ 

\begin{thm}
\label{thm:asymp expansion}Let $L$ be an ample line bundle and let
$K_{k}$ be the Bergman kernel of the Hilbert space $\mathcal{H}(X,L^{k}).$
Any interior point in $D\cap X(0)$ has a neighbourhood where $K_{k}(x,y)e^{-k\phi(x)/2}e^{-k\phi(y)/2}$
admits an asymptotic expansion as\begin{equation}
k^{n}(\det(dd^{c}\phi)(x)+b_{1}(x,y)k^{-1}+b_{2}(x,y)k^{-2}+...)e^{k\phi(x,y)},\label{eq:exp in prop}\end{equation}
where $b_{i}$ are global well-defined functions expressed as polynomials
in the covariant derivatives of $dd^{c}\phi$ (and of the curvature
of the metric $\omega$) which can be obtained by the recursion given
in \cite{b-b-s}.
\end{thm}
\begin{proof}
The proof is obtained by adapting the construction in \cite{b-b-s},
concerning \emph{positive} Hermitian line bundles, to the present
situation. The approach in \cite{b-b-s} is to first construct a {}``local
asymptotic Bergman kernel'' close to any point where $\phi$ is smooth
and $dd^{c}\phi>0$. Hence, the local construction applies to the
present situation as well. Then the local kernel is shown to differ
from the true kernel by a term of order $O(k^{-\infty}),$ by solving
a $\overline{\partial}$-equation with a good $L^{2}-$estimate. This
is possible since $dd^{c}\phi>1/C$ \emph{globally} in that case.
In the present situation we are done if we can solve \begin{equation}
\overline{\partial}u_{k}=g_{k},\label{eq:inhom dbar}\end{equation}
 where $g_{k}$ is a $\overline{\partial}-$closed $(0,1)-$form with
values in $L^{k},$ supported on the interior of the bounded set $D\cap X(0),$
with an estimate \begin{equation}
\left\Vert u_{k}\right\Vert _{k\phi}\leq C\left\Vert g_{k}\right\Vert _{k\phi}\label{eq:horm est}\end{equation}
To this end we apply the $L^{2}-$estimates of Hörmander-Kodaira \cite{de4}
with the metric $\psi_{k}$ on $L^{k}$ (formula \ref{pf of thm lnk: metric})
occurring in the proof of theorem \ref{thm:ln K}. This gives a solution
$u_{k}$ satisfying the inequality \ref{eq:horm est} with $k\phi$
replaced by the weight function $\psi_{k}.$ To see that we actually
have the estimate \ref{eq:horm est} (i.e. with the weight $k\phi$
itself) we apply the same argument as in the end of the proof of theorem
\ref{thm:ln K} to the left hand side in \ref{eq:horm est}. Finally,
for the right hand side in \ref{eq:horm est}, we use that $\phi=\phi_{e}$
on the bounded set $D\cap X(0)$ where $g_{k}$ is supported. 
\end{proof}

\section{\label{sec:Sections-vanishing-along}Sections vanishing along a divisor}

In this section we will show that the setup in the previous sections
can be adapted to the case when the Hilbert space $\mathcal{H}(X,L^{k})$
is replaced by the subspace of all sections vanishing to order at
least $k$ along a fixed divisor (which we for simplicity take to
be irreducible): \[
Z=\{ s_{Z}=0\},\]
assuming that the line bundle $L\otimes[Z]^{-1}$ is ample (without
assuming that $L$ is ample). The point is that this amounts essentially
to replacing $(L,\phi)$ with the Hermitian line bundle $(L\otimes[Z]^{-1},\phi-\ln(\left|s_{Z}\right|^{2}))$
to which the previous setup essentially applies if one also takes
into account the fact the singular metric $\phi-\ln(\left|s_{Z}\right|^{2})$
is equal to infinity on $Z.$ Hence, the corresponding curvature current
will be \emph{negative} close to $Z,$ even if the curvature form
of $\phi$ is globally positive.

\subsection{Equilibrium metrics with poles along a divisor}

Let $\mathcal{L}_{(X,L);Z}$ be the subclass of $\mathcal{L}_{(X,L)}$
consisting of all metrics $\widetilde{\phi}$ on $L$ such that the
Lelong numbers of $\widetilde{\phi}$ are bounded from below by one
along the divisor $Z:$\[
\nu(\widetilde{\phi})_{x}\geq1,\,\,\,\textrm{when\,}x\in Z,\]
 where\[
\nu(\widetilde{\phi})_{x}:=\lim_{r\rightarrow0_{+}}\frac{1}{r^{2n}}\int_{\left|z\right|\leq r}dd^{c}\phi\wedge(dd^{c}\left|z\right|^{2})^{n-1}/(n-1)!,\]
 with respect to any local coordinate $z$ centered at $x.$ Then
we define the associated \emph{equilibrium metric with poles along
$Z$} as \[
\phi_{e.Z}(x)=\sup\left\{ \widetilde{\phi}(x):\,\widetilde{\phi}\in\mathcal{L}_{(X,L);Z},\,\widetilde{\phi}\leq\phi\,\,\textrm{on$\, X$}\right\} \]
and the following set, compactly included in $X-Z:$\begin{equation}
D_{Z}:=\left\{ \phi_{e,Z}=\phi\right\} \label{eq:def of DZ}\end{equation}

\begin{lem}
\label{lem:The-following-decomposition}The following decomposition
holds: \[
\phi_{e,Z}=\psi_{e}+\ln(\left|s_{Z}\right|^{2})\]
where \[
\psi_{e}(x)=\sup\left\{ \widetilde{\psi}(x):\,\widetilde{\psi}\in\mathcal{L}_{(X,L\otimes[Z]^{-1})},\,\widetilde{\psi}\leq\phi-\ln(\left|s_{Z}\right|^{2})\,\,\textrm{on$\, X-Z$}\right\} .\]

\end{lem}
\begin{proof}
First observe that \begin{equation}
\nu(\widetilde{\phi})_{x}\geq1,\,\forall x\in Z\Leftrightarrow\exists(C,U_{Z}):\,\widetilde{\phi}\leq\ln(\left|s_{Z}\right|^{2})+C\,\textrm{on\,$U_{Z},$ }\label{eq:pf of lemma lelong 1}\end{equation}
where $U_{Z}$ is a neighbourhood of $Z.$ Indeed, this is a direct
consequence of the following characterization \cite{de4} of the Lelong
number at $0$ of a germ of a psh function in $\C^{n}:$ \begin{equation}
\nu(\widetilde{\phi})_{0}=\sup\left\{ \gamma:\,\widetilde{\phi}\leq\gamma\ln(\left|z\right|^{2})+C_{\gamma}\,\textrm{close to $0$}\right\} \label{eq:char of lelong}\end{equation}
Next, observe that \begin{equation}
\exists(C,U_{Z}):\,\widetilde{\phi}\leq\ln(\left|s_{Z}\right|^{2})+C\,\textrm{on\,$U_{Z}\Leftrightarrow\widetilde{\psi}:=\widetilde{\phi}-\ln(\left|s_{Z}\right|^{2})\in\mathcal{L}_{(X,L\otimes[Z]^{-1})}.$ }\label{eq:pf of lemma lelong 2}\end{equation}
To see this first note that in general $\widetilde{\psi}$ above determines
a (possibly singular) metric on $L\otimes[Z]^{-1}$ over $X-Z$ with
positive curvature current $dd^{c}\widetilde{\psi}=dd^{c}\widetilde{\phi}.$
Now the condition that $\widetilde{\phi}\leq\ln(\left|s_{Z}\right|^{2})+C$
close to $Z$ means that $\widetilde{\psi}$ is bounded close to $Z,$
which is equivalent to the fact that it extends to an element of $\mathcal{L}_{(X,L\otimes[Z]^{-1})}.$
This follows from the corresponding local extension property of plurisubharmonic
functions over an analytic variety (or more generally over any pluripolar
set) in $\C^{n}$ \cite{kl}.

Finally, combining \ref{eq:pf of lemma lelong 1} and \ref{eq:pf of lemma lelong 2}
gives a bijection between $\mathcal{L}_{(X,L);Z}$ and $\mathcal{L}_{(X,L\otimes[Z]^{-1})}.$
Since, $\widetilde{\phi}\leq\phi$ on $X$ if and only if the equality
holds on $X-Z$ (using that $\widetilde{\phi}$ is assumed to be u.s.c)
this finishes the proof of the lemma.
\end{proof}
Using the decomposition in the previous lemma the following regularity
result is obtained (compare theorem \ref{thm:reg}):

\begin{thm}
\label{thm:ref div}Suppose that $L\otimes[Z]^{-1}$ is an ample line
bundle and that the given metric $\phi$ on $L$ is smooth. Then 

(a) $\phi_{e,Z}$ is locally in the class $\mathcal{C}^{1,1}$ on
$X-Z,$ i.e. $\phi_{e}$ is differentiable and all of its first partial
derivatives are locally Lipschitz continuous. In fact, \[
\phi_{e,Z}=\psi_{e}+\ln(\left|s_{Z}\right|^{2}),\]
 where $\psi_{e}$ is locally in the class $\mathcal{C}^{1,1}$ on
all of $X.$

(b) The Monge-Ampere measure of $\phi_{e}$ on $X-Z$ is absolutely
continuous with respect to any given volume form and coincides with
the corresponding $L_{loc}^{\infty}$ $(n,n)-$form obtained by a
point-wise calculation: \begin{equation}
(dd^{c}\phi_{e,Z})^{n}=\det(dd^{c}\phi_{e,Z})\omega_{n}\,\,\,\textrm{on\,$X-Z$}\label{eq:ptwise repr of equil meas}\end{equation}

(c) the following identity holds almost everywhere on the set $D_{Z}=\{\phi_{e,Z}=\phi\}:$
\begin{equation}
\det(dd^{c}\phi)=\det(dd^{c}\phi_{e,Z})=\det(dd^{c}\psi_{e})\label{eq:monge on D}\end{equation}

\end{thm}
\begin{proof}
(a) By the previous lemma it is enough to prove that $\psi_{e}$ is
locally in the class $\mathcal{C}^{1,1}$ on $X.$ Since $L\otimes[Z]^{-1}$
is an ample line bundle over $X$ this would be a direct consequence
of theorem \ref{thm:reg} if $\psi:=\phi-\ln(\left|s_{Z}\right|^{2})$
were in the class $\mathcal{C}^{2}$ on all of $X$ and not only on
$X-Z.$ In order to supply the necessary modifications of the argument
note that, as a locally psh function, $\psi_{e}$ is necessarily bounded
from above close to $Z$ and hence the bound corresponding to the
bound \ref{pf thm reg: bound on cand} trivially holds over a neighbourhood
of $Z,$ since $\psi$ is equal to infinity along $Z.$ This proves
that $\psi_{e}$ is locally Lipschitz and the rest of the proof of
(a) can be modified an a similar way. The proof of $(b)$ and $(c)$
proceeds exactly as in the proof of theorem \ref{thm:reg}.
\end{proof}

\subsection{Bergman kernels vanishing along a divisor}

Now consider the sub Hilbert space $\mathcal{H}_{k,Z}$ of $\mathcal{H}(X,L^{k}),$
consisting of all sections $\alpha_{k}$ such that the order of vanishing
of $\alpha_{k}$ along the divisor $Z$ is at least $k.$ Since the
later condition means that \[
f_{k}=g_{k}\otimes s_{Z}^{\otimes k},\,\,\, g_{k}\in H^{0}(X,L^{k}\otimes[Z]^{-1})\]
the Hilbert space $\mathcal{H}_{k,Z}$ is isomorphic to the vector
space $H^{0}(X,L^{k}\otimes[Z]^{-1})$ equipped with the norm \ref{eq:norm restr}
induced by $\phi,$ under the natural embedding \begin{equation}
(\cdot)\otimes s_{Z}^{\otimes k}:\,\, H^{0}(X,L^{k}\otimes[Z]^{-1})\rightarrow H^{0}(X,L^{k})\label{eq:embedding of hilb}\end{equation}
We denote by $K_{k,Z}$ and $B_{k,Z}$ the corresponding Bergman kernels
and Bergman functions, respectively.

\begin{thm}
Assume that the line bundle $L\otimes[Z]^{-1}$ over $X$ is ample.
Let $B_{k,Z}$ be the Bergman function of the Hilbert space $\mathcal{H}_{k,Z}$
of all sections vanishing along $Z$ to order at least $k.$ Then\begin{equation}
k^{-n}B_{k,Z}\rightarrow1_{D_{e,Z}\cap X(0)}\det(dd^{c}\phi),\label{eq:l1 conv of B}\end{equation}
 in $L^{1}(X,\omega_{n}).$
\end{thm}
\begin{proof}
First note that the local Morse inequalities \ref{lem:(Local-Morse-inequalities)}
still hold when $B_{k}$ is replaced by $B_{k,Z},$ since $B_{k,Z}\leq B_{k}.$
Moreover, if $f_{k}$ is a local representation of en element of $\mathcal{H}_{k,Z},$
the embedding \ref{eq:embedding of hilb} and the characterization
\ref{eq:char of lelong} of the Lelong numbers give that $\frac{1}{k}\ln\left|f_{k}(z)\right|^{2}$(and
hence also $k^{-1}\textrm{ln\,$K_{k,Z}$) }$belongs to the class $\mathcal{L}_{(X,L);Z}.$
Hence, the proof of lemma \ref{eq:exp decay} goes through in the
present setting, showing that $k^{-n}B_{k,Z}$ converges (exponentially)
to zero on $X-D_{Z}.$ Thus, we get as in the proof of theorem \ref{thm:B in L1}
\begin{equation}
\lim_{k}k^{-n}\dim\mathcal{H}_{k,Z}=\lim_{k}\int_{D_{Z}}k^{-n}B_{k,Z}\omega_{n}\leq\int_{D_{Z}\cap X(0)}(dd^{c}\phi)^{n}/n!\label{pf of thm B div}\end{equation}
By theorem \ref{thm:ref div} the right hand side is bounded by $\int_{X-Z}(dd^{c}\psi_{e})^{n}/n!$
which is equal to the top intersection number of $c_{1}(L\otimes[Z]^{-1}).$
Now, since $L\otimes[Z]^{-1}$ is ample this number is equal to the
limit of $k^{-n}\dim H^{0}(X,L^{k}\otimes[Z]^{-1}),$ which in turn
equals the left hand side in \ref{pf of thm B div}. Hence, the inequality
in \ref{pf of thm B div} is actually an equality. The rest of the
proof proceeds word for word as in the proof of theorem \ref{thm:B in L1}.
\end{proof}
Next, we have the following generalization of theorem \ref{thm:ln K},
which can be seen as a global version of the general approximation
results for $(1,1)-$currents of Demailly \cite{de2}, in the particular
case when the current (here given by $dd^{c}\phi_{e,Z})$ is singular
along a divisor.

\begin{thm}
Assume that the line bundle $L\otimes[Z]^{-1}$ over $X$ is ample.
Let $K_{k,Z}$ be the Bergman kernel of the Hilbert space $\mathcal{H}_{k,Z}$
of all sections vanishing along $Z$ to order at least $k.$ Then
\end{thm}
\[
k^{-1}\textrm{ln\,$K_{k,Z}(x,x)\rightarrow\phi_{e,Z}(x)$}\]
uniformly on $X-Z$ (the rate of convergence is of the order $\ln k/k).$
In particular, the corresponding convergence of the Lelong numbers
along $Z$ holds: \[
\nu(k^{-1}\textrm{ln\,$K_{k,Z})_{x}\rightarrow\nu(\phi_{e,Z})_{x}(=1)$}\]
and the corresponding {}``Bergman volume forms'' converge on $X-Z$
to the equilibrium measure: 

\[
(dd^{c}(k^{-1}\textrm{ln\,$K_{k,Z}(x,x)))^{n}\rightarrow(dd^{c}\phi_{e,Z})^{n}$}\]
weakly as measures on $X-Z.$

\begin{proof}
As above the proof can be reduced to the proof of the corresponding
theorem in the case when there is no divisor $Z$ (theorem \ref{thm:ln K}).
This gives the uniform convergence \[
k^{-1}\textrm{ln\,$K_{k,Z}(x,x)-\ln(\left|s_{Z}\right|^{2}(x))\rightarrow\phi_{e}(x)-\ln(\left|s_{Z}\right|^{2}(x))$}\]
 on all of $X.$ The convergence of the Lelong numbers then follows
from the characterization \ref{eq:char of lelong}.
\end{proof}
Finally, note that the theorems \ref{thm:k as meas} and \ref{thm:asymp expansion}
generalize in the corresponding way to the present situation. The
explicit statements are omitted.

\section{\label{sec:Examples}Examples}

Finally, we illustrate some of the previous results with the following
examples, which can be seen as variants of the setting considered
in \cite{berm4} (compare remark \ref{rem:The-setting-considered}
below).

\begin{example}
\label{exa:Pn}Let $X$ be the $n-$dimensional projective space $\P^{n}$
and let $L$ be the hyperplane line bundle $\mathcal{O}(1)$. Then
$H^{0}(X,L^{k})$ is the space of homogeneous polynomials in $n+1$
homogeneous coordinates $Z_{0},Z_{1},..Z_{n}.$ The Fubini-Study metric
$\phi_{FS}$ on $\mathcal{O}(1)$ may be suggestively written as $\phi_{FS}(Z)=\ln(\left|Z\right|^{2})$
and the Fubini-Study metric $\omega_{FS}$ on $\P^{n}$ is the normalized
curvature form $dd^{c}\phi_{FS}.$ Hence the induced norm on $H^{0}(X,L^{k})$
is invariant under the standard action of $SU(n+1)$ on $\P^{n}.$
We may identify $\C^{n}$ with the {}``affine piece'' $\P^{n}-H_{\infty}$
where $H_{\infty}$ is the {}``hyperplane at infinity'' in $\C^{n}$
(defined as the set where $Z_{0}=0).$ In terms of the standard trivialization
of $\mathcal{O}(1)$ over $\C^{n}$ (obtained by setting $Z_{0}=1)$
the space $H^{0}(Y,L^{k})$ may be identified with the space of polynomials
$f_{k}(\zeta)$ in $\C_{\zeta}^{n}$ of total degree at most $k$
and the metric $\phi_{FS}$ on $\mathcal{O}(1)$ may be represented
by the function\[
\phi_{FS}(\zeta)=\ln(1+\left|\zeta\right|^{2}).\]
Moreover, any smooth metric on $\mathcal{O}(1)$ may be represented
by a function $\phi(\zeta)$ satisfying the following necessary growth
condition%
\footnote{in order that $\phi$ extend over the hyperplane at infinity to a
\emph{smooth} metric further conditions are needed.%
} \begin{equation}
-C+\ln(1+\left|\zeta\right|^{2})\leq\phi(\zeta)\leq\ln(1+\left|\zeta\right|^{2})+C,\label{eq:growth cond}\end{equation}
 which makes sure that the norm \ref{eq:norm restr}, expressed as
\[
\left\Vert f_{k}\right\Vert _{k\phi}^{2}:=\int_{\C^{n}}\left|f_{k}(\zeta)\right|^{2}e^{-k\phi(\zeta)}\omega_{FS}^{n}/n!\]
 is finite precisely when $f_{k}$ corresponds to a section of $\mathcal{O}(m)$,
for $m=1.$ In particular, any smooth compactly supported function
$\chi(\zeta)$ determines a \emph{smooth} perturbation \begin{equation}
\phi_{\chi}(\zeta):=\phi_{FS}(\zeta)+\chi(\zeta)\label{eq:pert metric}\end{equation}
of $\phi_{FS}$ on $\mathcal{O}(1)$ over $\P^{n},$ to which the
results in section \ref{sec:Equilibrium-measures-for} and \ref{sec:Bergman-kernels}
apply. For example, if $\chi$ is a radial function then it can be
checked that the graph of the equilibrium metric determined by $\phi_{\chi}$
is simply the convex hull of the graph of $\phi_{\chi}$ considered
as a function of $v:=\ln\left|\zeta\right|^{2}.$
\end{example}
The next example introduces a divisor into the picture, as in section
\ref{sec:Sections-vanishing-along}.

\begin{example}
\label{exa:div at inf}Let $X=\P^{n},$ $L=\mathcal{O}(2)$ and denote
by $Z$ the hyperplane at infinity in $\C^{n}.$ Then $L\otimes[Z]^{-1}\simeq\mathcal{O}(1)$
is ample. We equip $L$ with the canonical metric $2\phi_{FS}.$ Then
$H^{0}(X,L^{k})$ may be identified with the space of all polynomials
in $\C^{n}$ of degree at most $2k,$ while the subspace $\mathcal{H}_{k,Z}$
is the space of polynomials of degree at most $k.$ In this case the
set $D_{Z}$ (formula \ref{eq:def of DZ}) is, by symmetry, a ball
$B(0;r)$ centered at $0$ of radius $r,$ where $r$ is determined
by the following volume condition: \[
\int_{B(0;r)}(2\omega_{FS})^{n}/n!=1.\]

\end{example}
\begin{rem}
\label{rem:The-setting-considered}The setting considered in \cite{berm4}
corresponds to replacing $\omega_{FS}$ by the Euclidean metric on
$\C^{n}$ and the growth-condition \ref{eq:growth cond} by the condition
that \begin{equation}
(1+\epsilon)\ln(1+\left|\zeta\right|^{2})\leq\phi(\zeta)\label{eq:growth with eps}\end{equation}
(but with no upper bound or further assumption about smoothness {}``at
infinity''). Hence, in a certain sense, example \ref{exa:Pn} may
be seen as a limiting case of the setting in \cite{berm4}. In particular,
the set $D\cap\C^{n}$ may be non-compact in that example, while it
is always compact under the assumption \ref{eq:growth with eps}.
Note however, that example \ref{exa:div at inf} is essentially a
special case of the situation studied in \cite{berm4} (if $\omega_{FS}$
is replaced by the Euclidean metric). Indeed, the bound \ref{eq:growth with eps}
holds with $\epsilon=2.$ 
\end{rem}
In the next example, which in a certain sense is dual to the previous
one, the line bundle $L$ is not ample, but since $L\otimes[Z]^{-1}$
is, the main results in section \ref{sec:Sections-vanishing-along}
still apply.

\begin{example}
Let $X=\widetilde{\P^{n}}$ be the blow-up of $\P^{n}$ at the origin
in $\C^{n}$ and denote by $\pi$ the projection (blow-down map) from
$\widetilde{\P^{n}}$ to $\P^{n}.$ Let $L=\pi^{*}\mathcal{O}(2)$
and denote by $Z$ the exceptional divisor over $0.$ Then $L\otimes[Z]^{-1}$
is ample (for example by the Nakai-Moishezon criteria \cite{gr-ha}).
We equip $L$ with the metric $\pi^{*}(2\phi_{FS}).$ Then $H^{0}(X,L^{k})$
may again be identified with the space of all polynomials in $\C^{n}$
of total degree at most $2k,$ while the subspace $\mathcal{H}_{k,Z}$
is the space of polynomials of total degree at least $k+1.$ In this
case the set $D_{Z}$ (formula \ref{eq:def of DZ}) is, by symmetry,
the \emph{complement} of a ball $B(0;r)$ of radius $r,$ where $r$
is determined by the following volume condition: \[
\int_{\C^{n}-B(0;r)}(2\omega_{FS})^{n}/n!=1.\]
In order to compare with example \ref{exa:Pn}, where there is no
divisor, note that the blow-down map $\pi$ induces an isomorphism
$H^{0}(X,(L\otimes[Z]^{-1})^{k})\simeq H^{0}(\P^{n},\mathcal{O}(1)^{k}$
such that the {}``push-forward'' of the corresponding singular metric
$\psi$ on $L\otimes[Z]^{-1}$ (compare the decomposition in lemma
\ref{lem:The-following-decomposition}) becomes the singular metric
\[
2\phi_{FS}(\zeta)-\ln\left|\zeta\right|^{2}\]
on $\mathcal{O}(1).$ This metric may be seen as a limit of metrics
$\phi_{\chi_{i}}$ on $\mathcal{O}(1)$ of the form \ref{eq:pert metric},
where the limiting function $\chi$ is given by \[
\chi(\zeta)=\ln(1+\left|\zeta\right|^{-2}).\]
 The point is that $\chi$ tends to $0$ as $\left|\zeta\right|$
tends to infinity and to infinity as $\left|\zeta\right|$ tends to
$0.$
\end{example}

\end{document}